\documentclass[12pt]{amsart}
\usepackage{amsfonts,amssymb,amsthm,eucal,amsmath}
\usepackage{xy}
\input xy
\xyoption{all}
\allowdisplaybreaks
\newtheorem{thm}{Theorem}
\newtheorem{cor}[thm]{Corollary}
\newtheorem{lemma}[thm]{Lemma}

\newcommand{\R}{\mathbb{R}}
\newcommand{\C}{\mathbb{C}}

\newcommand{\E}{\mathbb{E}}
\newcommand{\N}{\mathbb{N}}

\newcommand{\ds}{\displaystyle}

\renewcommand{\P}{\mathbb{P}}

\renewcommand{\L}{\mathcal{L}}
\newcommand{\I}{\mathbb{I}}

\newcommand{\n}{\mathfrak{N}}
\renewcommand{\O}{\mathcal{O}}
\newcommand{\U}{\mathcal{U}}
\DeclareMathOperator{\tr}{Tr\,}

\renewcommand{\Re}{\operatorname{Re\,}}

\begin{document}

\title{Linear Functions on the Classical Matrix Groups}
\author[E.\ Meckes]{Elizabeth Meckes}\thanks{Research supported in
part by the ARCS Foundation. {\em AMS 2000 subject classifications.}
Primary 60F05; secondary 60B15, 60B10.  
}

\maketitle

\begin{quote}
\begin{center}{\sc Abstract}\end{center}
{\small Let $M$ be a random matrix in the orthogonal group $\O_n$, distributed
according to Haar measure, and let $A$ be a fixed $n\times n$ matrix
over $\R$ such that $\tr(AA^t)=n$.  Then the total variation distance of
the random variable $\tr(AM)$ to a standard normal random variable 
is bounded by $\frac{2\sqrt{3}}
{n-1}$, and this rate is sharp up to the constant.  Analogous 
results are obtained for $M$ a random unitary 
matrix and $A$ a fixed $n\times n$ matrix over $\C$.  
The proofs are applications of a new abstract normal approximation theorem 
which extends Stein's method of
 exchangeable pairs to situations in which continuous symmetries are present.}

\end{quote}

\medskip


\section{Introduction}
Let $\O_n$ denote the group of $n\times n$ orthogonal matrices, and let 
$M$ be distributed according to Haar measure on $\O_n$.  
Let $A$ be a fixed $n\times n$ matrix over $\R$, subject to the condition
that $\tr(AA^t)=n$, and let $W=\tr(AM)$.  D'Aristotile, 
Diaconis, and Newman showed in \cite{limit} that 
$$\sup_{\substack{\tr(AA^t)=n
\\-\infty<x<\infty}}\left|\P(W\le x)-\Phi(x)\right|\to0$$
as $n\to\infty$.  Their argument uses classical methods involving 
sub-subsequences and tightness, and cannot be improved to yield
a theorem for finite $n$.  
Theorem \ref{orth-main1} below gives an explicit rate of 
convergence of the law of $W$ to the standard normal distribution
in the total variation
metric on probability measures, specifically,
\begin{equation}\label{mainst}
d(\L_W,\n(0,1))_{TV}\le\frac{2\sqrt{3}}{n-1}\end{equation}
for all $n\ge2.$

The history of this problem begins with the following theorem, first 
given rigorous proof by Borel in
\cite{borel}: let $X$ be a random vector on the unit sphere $S^{n-1}$, and
let $X_1$ be the first coordinate of $X$.  Then $\P\big(\sqrt{n}X_1\le t\big)
\longrightarrow\Phi(t)$ as $n\to\infty,$ where $\Phi(t)=\frac{1}{\sqrt{2\pi}}
\int_{-\infty}^te^{-\frac{x^2}{2}}dx.$  Since the first column of a 
Haar-distributed orthogonal matrix is uniformly distributed on the unit 
sphere, Borel's theorem follows from Theorem \ref{orth-main1} by taking 
$A=\sqrt{n}\oplus{\bf0}$.
Borel's theorem was generalized in one direction by Diaconis and 
Freedman \cite{diafree}, who proved the convergence of the first $k$ 
coordinates of $\sqrt{n}X$ to independent standard normal random variables
 in total variation
distance for $k=o(n)$; \cite{diafree} also contains a detailed history 
of this problem.  This line of research was further developed in
\cite{del}, where a total variation bound was given between an
$r\times r$ block of a random orthogonal matrix and an $r\times r$ 
matrix of independent standard Gaussians, for $r=O\left(n^{1/3}\right).$
This was later improved by Jiang (see \cite{jia}) to 
$r=O\left(n^{1/2}\right)$, which he proved was sharp.  
In the same paper, Jiang also showed that given a sequence of Haar distributed random matrices $\{M_n\}$, there is a sequence of Gaussian matrices $\{Y_n\}$ with $Y_j$ defined on the same probability space as $M_j$ such that if 
$$\epsilon_n=\max_{\substack{1\le i\le n\\ 1\le j\le m_n}}\left|\sqrt{n}M_{ij}-Y_{ij}\right|$$
with $m_n\le\frac{n}{\log^2n},$
then $\epsilon_n\to0$ in probability as $n\to\infty$.  Thus an 
$n\times\frac{n}{\log^2n}$ block of a Haar distributed matrix can be approximated by a Gaussian matrix `in probability'.
Theorem \ref{orth-main1}
 gives another
sense in which a random orthogonal matrix is close to a matrix of
independent normals by giving a uniform bound of distance to 
normal over all linear combinations of entries of $M$.

Another special case of Theorem \ref{orth-main1} is $A=I$, so that $W=\tr(M)$.  
Diaconis and Mallows (see \cite{diamal}) first proved that $\tr(M)$ is 
approximately normal; Stein \cite{steintech} and Johansson \cite{joh} later 
independently obtained fast 
rates of
convergence to normal of $\tr(M^k)$ for fixed $k$, with Johansson's rates 
an improvement on Stein's.  In studying eigenvalues of random orthogonal 
matrices, Diaconis and Shahshahani \cite{diasha} extended this to show that
the joint limiting distribution of $\tr(M), \tr(M^2),\ldots, \tr(M^k)$ 
converges to that of independent normal variables as $n\to\infty$, for
$k$ fixed.  

The other source of motivation for theorems like Theorem \ref{orth-main1} is
Hoeffding's combinatorial central limit theorem \cite{hoeff}, which can
be stated as follows.  Let $A=(a_{ij})$ be a fixed $n\times n$ 
matrix over $\R$, normalized to have row and column sums equal to zero and 
$\frac{1}{n-1}\sum_{i,j}a_{ij}^2=1$.  Let $\pi$ be a random permutation in
$S_n,$ and let $W(\pi)=\sum_i a_{i\pi(i)}.$  Then under certain conditions
on $A$, $W$ is approximately normal.  Later, Bolthausen \cite{bolt}
proved an explicit rate of convergence via Stein's method.  Note that if 
$$M_{ij}=\begin{cases}1&\pi(j)=i\\0&\mbox{ otherwise}\end{cases}$$
then $W=\tr(AM)$, and so Hoeffding's theorem is really a theorem about the
distribution of linear
functions on the set of permutation matrices.

The unitary group is another source of many important applications; see, e.g. 
\cite{dia}.  In Section \ref{unit}, the 
random variable $\tr(AM)$ for $A$ a fixed matrix over $\C$ and $M$ a random
unitary matrix distributed according to Haar measure on $\U_n$ is 
considered.  
The main theorem of the section, Theorem \ref{2main} gives a 
 bound on the total variation distance of 
$Re\big[\tr(AM)\big]$ to standard normal
 analogous to that of Theorem \ref{orth-main1}; this can be viewed as theorem 
about real-linear functions on $\U_n$.  Corollary \ref{complex} 
shows that in the limit, the 
complex random variable $\tr(AM)$ is close to standard complex normal.
The methods used here cannot be used directly to prove the convergence of
$\tr(AM)$ to the standard complex normal; they work for
approximation of real-valued random variables only.  
A version of the present methods in a multivariate
context is forthcoming in
\cite{CM}, which includes a rate of convergence for Corollary \ref{complex}.

{\bf Notation and Conventions.}
The total variation distance $d_{TV}(\mu,\nu)$ between the measures $\mu$ and
$\nu$ on $\R$ is defined by
$$d_{TV}(\mu,\nu)=\sup_A\big|\mu(A)-\nu(A)|,$$
where the supremum is over measurable sets $A$.  This is equivalent to
$$d_{TV}(\mu,\nu)=\frac{1}{2}\sup_{f}\left|\int f(t)d\mu(t)-
\int f(t)d\nu(t)\right|,$$
where the supremum is taken over continuous functions which are bounded
by 1 and vanish at infinity; this is the definition used in what follows. 
The total variation distance between two random variables $X$ and $Y$ is
defined to be the total variation distance between their distributions:
$$d_{TV}(X,Y)=\sup_A\big|\P(X\in A)-\P(Y\in A)\big|=\frac{1}{2}\sup_{f}
\big|\E f(X)-\E f(Y)\big|.$$
We will use $\n(\mu,\sigma^2)$ to denote the normal distribution on $\R$ with 
mean $\mu$ and variance $\sigma^2$.  

\medskip

{\bf Acknowledgements.}
I would like to thank Persi Diaconis for sharing his many
insights with me.

\section{An abstract normal approximation theorem}\label{chap-abs1}
In this section, a general approach for normal approximation to random
variables with continuous symmetries is developed.  The ideas which give 
rise to Theorem \ref{abscont} below first appeared in Stein \cite{steintech},
where fast rates of convergence to Gaussian (as $n\to\infty$) 
were obtained for $\tr(M^k)$,
with $k\in\N$ fixed and $M$ a random $n\times n$ orthogonal
matrix.  
\begin{thm}\label{abscont}
Suppose that $(W,W_\epsilon)$ is a family of exchangeable pairs defined
on a common probability space with $\E W=0$ and $\E W^2=\sigma^2$.
Suppose that there are functions $\alpha$ and $\beta$ with
$$\E|\alpha(\sigma^{-1}W)|<\infty,\qquad\E|\beta(\sigma^{-1}W)|<\infty,$$
and a constant $\lambda$ such that 
\begin{enumerate}
\item $$\frac{1}{\epsilon^2}\E\left[W_\epsilon-W\big|
W\right]=-\lambda W+o(1)\alpha(W),$$ \label{lin-cond}
\item $$\frac{1}{\epsilon^2}\E\left[(W_\epsilon-W)^2\big|
W\right]=2\lambda\sigma^2 +E\sigma^2+o(1)\beta(W),$$ \label{quad-cond}
\item $$\frac{1}{\epsilon^2}\E\left|W_\epsilon-W\right|^3
=o(1),$$ \label{cube-cond}
\end{enumerate}
where $o(1)$ refers to the limit as $\epsilon\to0$, with 
the implied constants deterministic.

Then
$$d_{T.V.}(W,Z)\le \frac{1}{\lambda}\E\left|E\right|,$$
where $Z\sim\n(0,\sigma^2)$.
\end{thm}

\medskip

{\bf Remark:} The factor of $\frac{1}{\epsilon^2}$ in each of 
the three expressions above could be replaced by a general function 
$f(\epsilon)$.  In practice, $W_\epsilon$ is typically constructed such that
$W_\epsilon-W=O(\epsilon)$.  This makes it clear that
$f(\epsilon)=\frac{1}{\epsilon^2}$ is
the suitable choice for condition \ref{quad-cond}.  It is less clear that
$f(\epsilon)=\frac{1}{\epsilon^2}$ is the suitable choice for 
condition \ref{lin-cond}.  In the applications given here, while 
$W_\epsilon-W=O(\epsilon)$, symmetry conditions imply that
$$\E\left[W_\epsilon-W\big|W\right]=O(\epsilon^2).$$

\bigskip

Before beginning the proof, some background on Stein's method is helpful.
The following lemma is key.
\begin{lemma}[Stein]\label{gauss}Let $Z\sim\n(0,1).$ 
Then
\begin{enumerate}
\item For all $f\in C^1_o(\R),$
$$\E\big[f'(Z)-Zf(Z)\big]=0.$$
\item \label{converse}If $Y$ is a random variable such that 
$$\E\big[f'(Y)-Yf(Y)\big]=0$$
for all $f\in C_b^1(\R),$ then $\L(Y)=\L(Z);$ i.e., $Y$ is also distributed
as a standard Gaussian random variable.
\item For $g:\R\to\R$ with $\E g(Z)<\infty$ given, the function
\begin{equation}\label{U_odef}
U_og(t)=e^{t^2/2}\int_{-\infty}^t\big[g(x)-\E g(Z)\big]e^{-x^2/2}dx.
\end{equation}
is a solution to the differential equation
$$f'(x)-xf(x)=g(x)-\E g(Z).$$
\end{enumerate}
\end{lemma}

\medskip

The lemma says that the standard Gaussian distribution $\gamma$ on $\R$ is the 
unique distribution with the property that $\int_{\R}(f'(x)-xf(x))d\gamma(x)$ 
is always zero.  The idea of Stein's method is that
if $W$ is a random variable such that $\E\big[f'(W)-Wf(W)\big]$ is always
{\em small}, then the distribution of $W$ is close $\gamma$.  
There are several 
approaches to bounding this quantity; the approach taken here is modelled 
on the method of exchangeable pairs (see \cite{steinbook}).
In any of the approaches, the following bounds on $U_o$ are useful.

\begin{lemma}[Stein]\label{stbds}
Let $U_o$ be the operator defined in equation (\ref{U_odef}).  Then
\begin{enumerate}
\item $\|U_og\|_\infty\le \sqrt{\frac{\pi}{2}}\|g-\E g(Z)\|_\infty
\le \sqrt{2\pi}\|g\|_\infty$
\item $\|(U_og)'\|_\infty\le 2\|g-\E g(Z)\|_\infty\le 4\|g\|_\infty$
\item $\|(U_og)''\|_\infty\le 2\|g'\|_\infty$
\end{enumerate}
\end{lemma}

\medskip

With this background, the proof of Theorem \ref{abscont} is straightforward.
\begin{proof}[Proof of Theorem \ref{abscont}]
By considering $\sigma^{-1}W$ instead of $W$, we may without loss assume 
that $\sigma=1$.  
For $g\in C^\infty_o(\R)$ fixed, let $f$ be the solution given in equation 
\eqref{U_odef} to the differential 
equation
$$f'(x)-xf(x)=g(x)-\E g(Z).$$
Fix $\epsilon$.  By the exchangeability of $(W,W_\epsilon)$, 
\begin{equation}\label{exch}\begin{split}
0&=\E\left[(W_\epsilon-W)(f(W_\epsilon)+f(W))\right]\\
&=\E\left[(W_\epsilon-W)(f(W_\epsilon)-f(W))+2(W_\epsilon-W)f(W)\right]\\
&=\E\left[\E\left[(W_\epsilon-W)^2\big|W\right]f'(W)+2\E\left[(W_\epsilon-
W)\big|W\right]f(W)+R\right],
\end{split}\end{equation}
where $R$ is the error in the derivative approximation.  By Taylor's theorem
and Lemma \ref{stbds},
$$|R|\le\frac{\|f''\|_\infty}{2}|W_\epsilon-W|^3\le\|g'\|_\infty|W_\epsilon-W
|^3,$$ and so 
$$\lim_{\epsilon\to0}\frac{1}{\epsilon^2}\E|R|=0.$$
Dividing both sides of (\ref{exch}) by $2\lambda\epsilon^2$ and 
taking the limit 
as $\epsilon\to0$ gives:
$$0=\E\left[f'(W)-Wf(W)+\frac{E}{2\lambda}f'(W)\right]=\E\left[g(W)-g(Z)+
\frac{E}{2\lambda}f'(W)\right].$$
Rearranging and applying the bound on $\|f'\|$ 
from Lemma
\ref{stbds} yields
\begin{equation*}
\big|\E g(W)-\E g(Z)\big|\le\frac{2\|g\|_\infty}{\lambda}\E|E|.
\end{equation*}
Since $C^\infty_o(\R)$ is dense (with respect to the supremum norm)
in the class of bounded 
continuous functions vanishing at infinity, this completes the proof.
\end{proof}

\section{The Orthogonal Group}\label{orth}
This section is mainly devoted to the proof of the following theorem.

\begin{thm}\label{orth-main1}
Let $A$ be a fixed $n\times n$ matrix over $\R$ such that $\tr(AA^t)=n,$
 $M\in\O_n$ distributed according to Haar measure, and $W=\tr(AM)$.  Let $Z$
be a standard normal random variable. 
Then for $n>1$,
\begin{equation}\label{thm}
d(W,Z)_{TV}\le\frac{2\sqrt{3}}{n-1}.\end{equation}
\end{thm}
\vspace{5mm}
The bound in Theorem \ref{orth-main1} is sharp up to the constant; consider the
matrix $A=\sqrt{n}\oplus{\bf0}$ where ${\bf0}$ is the $n-1\times n-1$
matrix with all zeros.  For this $A$, Theorem \ref{orth-main1} reproves the following
theorem, proved in \cite{diafree} with slightly worse constant
\begin{thm}\label{sphere}
Let ${\bf x} \in\sqrt{n}S^{n-1}$ be uniformly distributed, and let $Z$ be
a standard normal random variable.  Then
$$d_{TV}(x_1,Z)\le\frac{2\sqrt{3}}{n-1}.$$
\end{thm}

It is shown in \cite{diafree} that the order of this error term is correct.

\begin{proof}[Proof of Theorem \ref{orth-main1}]
First note that one can
assume without loss of generality that $A$ is diagonal: let $A=UDV$ be the
singular value decomposition of $A$.  Then $W=\tr(UDVM)=\tr(DVMU)$, and the
distribution of $VMU$ is the same as the distribution of $M$ by the 
translation invariance of Haar measure.  

Now define the pair $(W,W_\epsilon)$ for
each $\epsilon$ as follows.  Choose $H=(h_{ij})\in\O(n)$ according to Haar 
measure, 
independent of $M$, and let $M_\epsilon=HA_\epsilon H^tM$, where
$$A_\epsilon=\begin{bmatrix}\sqrt{1-\epsilon^2}&\epsilon&&&\\-\epsilon&\sqrt
{1-\epsilon^2}&&0&\\&&1&&\\&0&&\ddots&\\&&&&1\end{bmatrix},$$
thus $M_\epsilon$ can be thought of as a small random rotation of $M$.  
Let $W_\epsilon=W(M_\epsilon)$; $(W,W_\epsilon)$ is an exchangeable pair
by construction.

It is convenient to rewrite $M_\epsilon$ as follows.  Let 
$I_2$ be the $2\times2$ identity matrix, $K$ the $n\times
2$ matrix consisting of the first two columns of $H$, and let 
$$C_2=\begin{bmatrix}0&1\\-1&0\end{bmatrix}.$$  Then 
\begin{equation*}\begin{split}
M_\epsilon&=M+K\big[(\sqrt{1-\epsilon^2}-1)I_2+\epsilon C_2\big]K^tM
\\&=M+K\left[\left(-\frac{\epsilon^2}{2}+O(\epsilon^4)\right)I_2+
\epsilon C_2\right]K^tM,\end{split}\end{equation*}
and so 
\begin{equation}
W_\epsilon-W=\epsilon\left[\left(-\frac{
\epsilon}{2}+O(\epsilon^3)\right)\tr(AKK^tM)+\tr(AKC_2K^tM)\right]
\label{diff}\end{equation}

Now, the distribution of $H$ is unchanged by multiplying a fixed 
row or column by $-1$ and $H$ is orthogonal, 
thus $\E h_{ij}h_{k\ell}=\frac{1}{n}\delta_{ik}
\delta_{j\ell}$.  This implies that 
$$\E \big[KK^t\big]=\frac{2}{n}I_n$$
and $$\E \big[KCK^t\big]=0;$$
combining this with
 (\ref{diff}) yields:
\begin{equation}\begin{split}
\frac{n}{\epsilon^2}\E&\left[(W_\epsilon-W)\big|W\right]\\&=
-\frac{n}{2}\E\left[\E\left[\tr(AKK^tM)\big|M\right]\big|W\right]+
\frac{n}{\epsilon}\E\left[\E\left[\tr(AKC_2K^tM)\big|M
\right]\big|W\right]+O(\epsilon)\\&=-\E\left[\E\left[\tr(AM)\big|M\right]
\big|W\right]+O(\epsilon)
\nonumber\\&=-W+O(\epsilon),
\end{split}\end{equation}
where the independence of $M$ and $H$ has been used to get the third line,
and the implied constants in the $O(\epsilon)$ here and in what follows  
may depend on $n$.
Condition \ref{lin-cond} of Theorem \ref{abscont} is thus 
satisfied with $\lambda=\frac{1}{n}$.

Recall now that $A$ is assumed to be diagonal.  
The second condition of Theorem \ref{abscont} can also be verified using
 the expression in (\ref{diff}) as follows.
\begin{equation}\begin{split}
\frac{n}{2\epsilon^2}&\E\left[(W_\epsilon-W)^2\big|W\right]\\
&=\frac{n}{2}\E\left[\E\left[(\tr(AKC_2K^tM))^2\big|M\right]\big|
W\right]+O(\epsilon)\\
&=\left.\frac{n}{2}\E\left[\sum_{i,j}\sum_{\substack{i'\neq i\\
j'\neq j}}m_{i'i}m_{j'j}a_{ii}a_{jj}\label{secdiff1}
\E\left[(h_{i1}h_{i'2}-h_{i2}h_{i'
1})(h_{j1}h_{j'2}-h_{j2}h_{j'1})\big|M\right]\right|W\right]+O(\epsilon),
\end{split}\end{equation} 
where the conditions on $i'$ and $j'$ are justified as the expression inside
the expectation is identically zero when either $i=i'$ or $j=j'$.  

Standard techniques are available for computing the mixed moments of 
entries of $H$; see e.g. \cite{hp}, section 4.2.  Using these techniques
 and the independence of $M$ and $H$ gives that
for $i'\neq i$ and $j'\neq j$,
\begin{equation}\label{useful-fourth}
\E\left[(h_{i1}h_{i'2}-h_{i2}h_{i'1})(h_{j1}h_{j'2}-h_{j2}h_{j'1})\big|M
\right]=\frac{2}{n(n-1)}\Big[\delta_{ij}\delta_{i'j'}-\delta_{ij'}
\delta_{ji'}\Big];\end{equation}
putting this into (\ref{secdiff1}) yields
\begin{eqnarray*}\nonumber
\frac{n}{2\epsilon^2}\E\left[(W_\epsilon-W)^2\big|W\right]
&=&\frac{1}{n-1}\sum_{i,j}\sum_{\substack{i'\neq i\\j'\neq j}}m_{i'i}m_{j'j}
a_{ii}a_{jj}\Big[\delta_{i'j'}\delta_{ij}-\delta_{ij'}\delta_{ji'}\Big]
+O(\epsilon)\\
\nonumber&=&\frac{1}{n-1}\left[\sum_ia_{ii}^2\left[(M^tM)_{ii}-m_{ii}^2\right]
-\sum_{i,i'\neq i}(MA)_{ii'}(MA)_{i'i}\right]+O(\epsilon)\\\nonumber&=&
\frac{1}{n-1}\left[n-\sum_ia_{ii}^2m_{ii}^2-\left[\tr((MA)^2)-\sum_ia_{ii}^2
m_{ii}^2\right]\right]+O(\epsilon)
\\&=&1+\frac{1}{n-1}
\left[1-\tr((AM)^2)\right]+O(\epsilon),
\end{eqnarray*}
thus 
\begin{equation}\label{secdiff2}
\lim_{\epsilon\to0}\frac{1}{\epsilon^2}\E\left[(W_\epsilon-W)^2\big|W\right]
=\frac{2}{n}+\frac{2}{n(n-1)}\left[1-\tr((AM)^2)\right]\end{equation}
and so 
\begin{equation}\label{orthE}E=\frac{2}{n(n-1)}\left[1-\tr((AM)^2)\right].
\end{equation}

Finally, \eqref{diff} gives immediately that 
$$\E\left[|W_\epsilon-W|^3\big|W\right]=O(\epsilon^3).$$
It remains to bound $n\E|E|.$

\begin{eqnarray}
\E\left[\tr((AM)^2)\right]&=&\E\left[\sum_{i,j}a_{ii}a_{jj}m_{ij}m_{ji}\right]
\nonumber\\
&=&\frac{1}{n}\sum_ia_{ii}^2=1\label{mean},
\end{eqnarray}
and
\begin{equation*}\begin{split}
\E\big[(\tr(&(AM)^2))^2\big]\\&=\E\left[\left(\sum_{i,j}a_{ii}a_{jj}m_{ij}
m_{ji}\right)\left(\sum_{k,l}a_{kk}a_{ll}m_{kl}m_{lk}\right)\right]\nonumber
\\\nonumber&=
\sum_{i,j,k,l}a_{ii}a_{jj}a_{kk}a_{ll}\left[\frac{n+1}{(n-1)n(n+2)}\left[
\delta_{ij}\delta_{kl}\big(1-\delta_{ik}\big)+
\delta_{ik}\delta_{jl}\big(1-\delta_{ij}\big)+\right.\right.
\\\nonumber&\left.\left.\qquad\qquad\qquad\qquad\qquad
\delta_{il}\delta_{jk}\big(1-\delta_{ij}\big)
\right]+\frac{3}{n(n+2)}\I(i=j=k=l)\right]\\
&=\frac{n+1}{n(n-1)(n+2)}\left(\sideset{}{'}\sum_{i,k}a_{ii}^2
a_{kk}^2+\sideset{}{'}\sum_{i,j}a_{ii}^2a_{jj}^2+\sideset{}{'}\sum_{
i,j}a_{ii}^2a_{jj}^2\right)+\frac{3}{n(n+2)}\sum_{i}a_{ii}^4\nonumber
\end{split}\end{equation*}
Now,
\begin{equation*}\begin{split}
\sideset{}{'}\sum_{i,j}a_{ii}^2a_{jj}^2&=\sum_ia_{ii}^2(n-a_{ii}^2)
=n^2-\sum_ia_{ii}^4.
\end{split}\end{equation*}
Applying this above gives
\begin{equation}\begin{split}
\E\big[(\tr((AM)^2))^2\big]&
=\frac{3(n+1)n^2}{(n-1)n(n+2)}-
\frac{3(n+1)}{(n-1)n(n+2)}\sum_ia_{ii}^4+\frac{3}{n(n+2)}\sum_ia_{ii}^4
\\&\le3+\frac{6}{(n-1)(n+2)}.\label{4th}
\end{split}\end{equation}

Putting these estimates into Theorem \ref{abscont} gives:
\begin{equation}d_{T.V.}(W,Z)\le\frac{2\sqrt{2+\frac{6}{(n-1)(n+2)}}}{(n-1)}.
\end{equation}
Noting that $\frac{6}{(n-1)(n+2)}\le1$ for
$n\ge3$ and that the bound in Theorem \ref{orth-main1} is trivially true for 
$n=2$ completes the proof.  
\end{proof}

\section{The Unitary Group}\label{unit}
Now let $M\in\U_n$ be distributed according 
to Haar measure, $A$ be an $n\times n$ matrix over $\C$, and
$W=\tr(AM)$.  In \cite{limit} it was shown that if $M=\Gamma+i\Lambda$ and 
$A$ and $B$ are fixed real diagonal matrices with $\tr (AA^t)=\tr(BB^t)=n,$
then $\tr(A\Gamma)+i\tr(B\Lambda)$ converges in distribution to a standard
complex normal random variable.  
This implies in particular that $\Re(W)$ converges in
distribution to $\n\left(0,\frac{1}{2}\right).$  The main theorem of 
this section gives a rate of this convergence in total 
variation distance.  

A more natural question might be the convergence of 
$W$ to a standard complex random variable.  As this is a multivariate problem,
Theorem \ref{abscont} cannot be applied.  A multivariate version of 
Theorem \ref{abscont}
is forthcoming in \cite{CM}, which also includes a rate of 
convergence of $W$ to 
a standard complex Gaussian random variable.

\begin{thm}\label{2main}
With $M$, $A$, and $W$ as above, let $W_\theta$ be the inner product of
$W$ with the unit vector making angle $\theta$ with the real axis.  Then
\begin{equation}
d_{TV}\left(W_\theta,\n\left(0,\frac{1}{2}\right)\right)\le\frac{c}{n}
\end{equation}
for a constant $c$ which is independent of $\theta$.
\end{thm}

The constant $c$ is asymptotically equal to $2\sqrt{2}$; for $n\ge 8$ it can
be taken to be 4.

\medskip

\begin{proof}
To prove the theorem, first note that it suffices to consider the case 
$\theta=0$, that is, to prove that 
$$d_{TV}\left(\Re(W),\n\left(0,\frac{1}{2}\right)\right)\le\frac{c}{n}.$$
The theorem then follows as stated since the distribution of $W$ is 
invariant under multiplication by any complex number of unit modulus.  Also,
$A$ can again be assumed diagonal with positive real entries by the singular 
value decomposition.

The proof is almost identical to the orthogonal case.  
Let $H\in\U_n$ be a random unitary matrix, independent of 
$M$, and let $M_\epsilon=HA_\epsilon H^*M$, where
$$A_\epsilon=\begin{bmatrix}\sqrt{1-\epsilon^2}&\epsilon&&&\\-\epsilon&\sqrt
{1-\epsilon^2}&&0&\\&&1&&\\&0&&\ddots&\\&&&&1\end{bmatrix}.$$
Let $W_\epsilon=W(M_
\epsilon)$.

Let
$I_2$ be the $2\times2$ identity matrix, $K$ the $n\times
2$ matrix consisting of the first two columns of $H$, and let 
$$C_2=\begin{bmatrix}0&1\\-1&0\end{bmatrix}.$$  Then
\begin{eqnarray}
W_\epsilon-W&=&\tr\left(\left(-\frac{\epsilon^2}{2}+O(\epsilon^4)\right)AKK^*M+
\epsilon AKC_2K^*M\right)\nonumber\\&=&\epsilon\left[\left(-\frac{
\epsilon}{2}+O(\epsilon^3)\right)\tr(AKK^*M)+\tr(AKC_2K^*M)\right].
\label{diff'}\end{eqnarray}

Let $W^r=Re(W)$ and $W^r_\epsilon=Re(W_\epsilon)$.  
As in the orthogonal case, to verify the conditions of Theorem \ref{abscont}
various mixed moments of the entries of $H$ are needed.  The relevant
unitary integrals can also be found in \cite{hp}, section 4.2.  They imply 
in particular that
\begin{eqnarray}
\E\left[KK^*_{ij}\right]&=&\frac{2}{n}\delta_{ij}\label{KK^*}\\
\E\left[KC_2K^*_{ij}\right]&=&0\label{KCK},
\end{eqnarray}
thus
\begin{equation}\label{firstdiff'}
\lim_{\epsilon\to0}\frac{n}{\epsilon^2}\E\left[W^r_\epsilon-W^r\big|W\right]=
-W^r;\end{equation}
condition \ref{lin-cond} is satisfied with $\lambda=\frac{1}{n}.$
Also by (\ref{diff'}),
\begin{equation}\begin{split}
\lim_{\epsilon\to0}\frac{n}{2\epsilon^2}\E[(W_\epsilon^r &
-W^r)^2\big|W]\\&=\lim_{\epsilon\to0}\frac{n}{2}\E\left[\left
(Re(\tr(AKC_2K^*M))\right)^2\big|W\right]\\
&=\frac{n}{4}Re\,\E\left[\sum_{i,j,k,l}a_{ii}m_{ji}a_{kk}m_{lk}(h_{i1}
\overline{h}_{j2}-h_{i2}\overline{h}_{j1})(h_{k1}
\overline{h}_{l2}-h_{k2}\overline{h}_{l1})\right.+
\\\label{secdiff1'}&\left.\left.\quad\qquad\qquad
\phantom{\sum_{i,j,k,l}}a_{ii}m_{ji}a_{kk}\overline{m}_{lk}(h_{i1}
\overline{h}_{j2}-h_{i2}\overline{h}_{j1})(\overline{h}_{k1}
h_{l2}-\overline{h}_{k2}h_{l1})\right|W\right]
\end{split}\end{equation}
Using the formulae from \cite{hp}, it is straightforward to show that
\begin{equation}\begin{split}\label{prod1'}
\E[(h_{i1}\overline{h}_{j2}-h_{i2}\overline{h}
_{j1})&(h_{k1}\overline{h}_{l2}-h_{k2}\overline{h}_{l1})]\\&=-\frac{2\delta_{il}
\delta_{jk}(1-\delta_{ij})}{(n-1)(n+1)}+
\frac{2\delta_{ij}\delta_{k\ell}(1-\delta_{ik})}{(n-1)n(n+1)}
-\frac{2\I(i=j=k=l)}{n(n+1)}
\end{split}\end{equation}
and
\begin{eqnarray}
\E\big[(h_{i1}\overline{h}_{j2}&&\hspace{-9mm}-\hspace{1mm}h_{i2}
\overline{h}_{j1})(\overline{h}_{k1}h_{l2}-\overline{h}_{k2}h_{l1})\big]
\hspace{3in}
\nonumber\\
&=&\frac{2\left(\delta_{ik}\delta_{jl}(1-\delta_{ij})
\right)}{(n-1)(n+1)}-\frac{2\left(\delta_{ij}\delta_{kl}(1-\delta_{ik})\right)}
{n(n-1)(n+1)}\label{prod2}+\frac{2\I(i=j=k=l)}{n(n+1)}.
\end{eqnarray}

Let $\ds\sideset{}{'}\sum_{i,j}$ stand for summing over all pairs $(i,j)$ where
$i$ and $j$ are distinct.  Putting (\ref{prod1'}) and (\ref{prod2}) into 
(\ref{secdiff1'}) and using the
independence of $M$ and $H$ gives:

\begin{equation*}\begin{split}
\lim_{\epsilon\to0}\frac{n}{2\epsilon^2}
\E\big[\big(&W_\epsilon^r-W^r\big)
^2\big|W\big]\\
&=\frac{n}{2(n-1)(n+1)}Re\,\E\left[\sum_{i,j,k,\ell}a_{ii}m_{ji}
a_{kk}m_{\ell k}\left(-\delta_{i\ell}\delta_{jk}(1-\delta_{ij})+
\frac{1}{n}\delta_{ij}\delta_{k\ell}(1-\delta_{ik})\right.\right.\\
&\hspace{3.5in}-\left.\left(\frac{n-1}{n}\right)\I(i=j=k=\ell)\right)
\\&\hspace{1.75in}
+\sum_{i,j,k,\ell}a_{ii}m_{ji}a_{kk}m_{\ell k}
\left(\delta_{ik}\delta_{j\ell}(1-\delta_{ij})-\frac{1}{n}\delta_{ij}
\delta_{k\ell}(1-\delta_{ik})\right.\\&\hspace{3.25in}
\left.\left.\phantom{\sum_i}\left.+\left(\frac{n-1}{n}\right)\I(i=j=k=l)
\right)
\right|W\right]
\\&=\frac{n}{2(n-1)(n+1)}Re\,\E\left[-\sideset{}{'}
\sum_{i,j}a_{ii}a_{jj}m_{ij}m_{ji}\right.+\frac{1}{n}
\sideset{}{'}\sum_{i,k}a_{ii}a_{kk}m_{ii}m_{kk}
\\&\hspace{2in}-\left(\frac{n-1}{n}\right)\sum_ia_{ii}^2m_{ii}^2+
\sideset{}{'}\sum_{i,j}a_{ii}^2|m_{ji}|^2
\\&\hspace{2in}-\frac{1}{n}\sideset{}{'}\sum_{i,k}a_{ii}a
_{kk}m_{ii}\overline{m}_{kk}\left.\left.+
\frac{n-1}{n}\sum_ia_{ii}^2|m_{ii}|^2\right|W\right]\\
&=\frac{n}{2(n-1)(n+1)}Re\,\E\left[-\left(\tr((AM)^2)-\sum_i(AM)_{ii}
^2\right)\right.+\frac{1}{n}\left(W^2-
\sum_i(AM)_{ii}^2\right)\\&\hspace{2in}
-\left(\frac{n-1}{n}\right)
\sum_i(AM)_{ii}^2+
\sum_ia_{ii}^2(1-|m_{ii}|^2)\\&\hspace{2in}\left.\left.
-\frac{1}{n}\left(\left|W\right|^2-\sum_ia_{ii}^2|m_{ii}|^2\right)
+\frac{n-1}{n}\sum_i
a_{ii}^2|m_{ii}|^2\right|W\right]
\\&=
\frac{1}{2}+\frac{1}{2(n-1)(n+1)}\\&
\left.+\frac{n}{2(n-1)(n+1)}Re\,\E\left[-\tr((AM)^2)
+\frac{W^2-|W|^2}{n}\right|W\right].
\end{split}\end{equation*}
Condition (2) of Theorem \ref{abscont} is thus satisfied with 
\begin{equation}\begin{split}\label{unitE}
nE&=\frac{1}{2(n-1)(n+1)}+\left.\frac{n}{2(n-1)(n+1)}
Re\,\E\left[-\tr((AM)^2)
+\frac{W^2-|W|^2}{n}\right|W\right].
\end{split}\end{equation}
It remains to estimate $n\E|E|$.  First,
\begin{eqnarray*}
\E\left|\tr((AM)^2)\right|&=&\E\sqrt{\sum_{i,j,k,l}a_{ii}a_{jj}m_{ij}m_{ji}
a_{kk}a_{ll}\overline{m}_{kl}\overline{m}_{lk}}\\
&\le&\sqrt{\sum_{i,j,k,l}a_{ii}a_{jj}a_{kk}a_{ll}
\E\left[m_{ij}m_{ji}\overline{m}_{kl}\overline{m}_{lk}\right]}\\
&=&\sqrt{\frac{2n^2}{(n-1)(n+1)}-\frac{2}{(n-1)n(n+1)}\left(\sum_i
a_{ii}^4\right)}\\&\le&\sqrt{2+\frac{1}{n^2-1}},
\end{eqnarray*}
using the formulae of \cite{hp} to evaluate the integrals.

Next,
\begin{eqnarray*}
\E|W|^2&=&\E\left[\sum_{i,j}a_{ii}a_{jj}m_{ii}\overline{m}_{jj}
\right]\\&=&\frac{1}{n}\sum_ia_{ii}^2\\&=&1.
\end{eqnarray*}

Putting these estimates into (\ref{unitE}) proves the theorem.

\end{proof}

\medskip

Theorem \ref{2main} yields the following bivariate corollary, which 
can also be seen as a corollary of the main unitary lemma of \cite{limit}.

\begin{cor}\label{complex}
Let $M$ be a random unitary matrix, $A$ a fixed $n\times n$ matrix over
$\C$ with $\tr(AA^*)=n$, and let $W=\tr(AM)$.  Then the distribution
of $W$ converges to the standard complex normal distribution in the weak-star
topology.
\end{cor}

\begin{proof}
The result follows immediately from Theorem \ref{2main} by considering 
the characteristic function of $W$.
\end{proof}

\bibliographystyle{plain}
\bibliography{matrices2}

\bigskip

\noindent Elizabeth S. Meckes\\
Department of Mathematics\\
Stanford University\\
Stanford, CA 94305\\

\noindent meckes@math.stanford.edu\\

\end{document}